\documentclass[english, 10pt]{article}

\usepackage{std_english}

\declaretheorem[title=Theorem, parent=section]{thm}

\declaretheorem[title=Proposition, numberlike=thm]{prop}

\declaretheorem[title=Lemma, numberlike=thm]{lem}

\declaretheorem[title=Question, numberlike=thm]{quest}

\declaretheorem[title=Conjecture, numberlike=thm]{conj}

\theorembodyfont{\rmfamily}

\declaretheorem[title=Proof,
				postheadhook=--~, 
				numbered=no,
				postfoothook=\hfill$\blacksquare$]{dem}

\theoremheaderfont{\mdseries\itshape}

\declaretheorem[title=Case,
				numberwithin=dem,
				]{case}

\declaretheorem[title=Claim,
				numberwithin=dem
				]{claim}

\newcommand{\ftwo}{\mathbb{F}_2}

\newcommand{\vts}[1]{V(#1)}
\newcommand{\eds}[1]{E(#1)}

\newcommand{\vn}{\varnothing}

\newcommand{\se}{\subseteq}

\newcommand{\sm}{\setminus}

\newcommand{\cro}[1]{\left[#1\right]}

\newcommand{\ceil}[1]{\lceil #1\rceil}

\newcommand{\rset}{\mathbb{R}}

\newcommand{\rsetp}{\mathbb{R}_+}

\newcommand{\set}[1]{\{#1\}}
\newcommand{\myref}[1]{\Cref{#1}\xspace}

\newcommand{\stab}[1]{\mathsf{STAB}(#1)}

\newcommand{\match}[1]{\mathsf{MATCH}(#1)}

\newcommand{\card}[1]{|#1|}

\newcommand{\one}{\mathbf{1}}

\newcommand{\idp}{integer decomposition property\xspace}

\newcommand{\cfp}{odd-$C_5^{+}$\xspace}
\newcommand{\cfpfree}{odd-$C_5^{+}$-free\xspace}

\newcommand{\np}{{\sf NP}\xspace}

\newcommand{\cspace}[1]{\mathcal{C}(#1)}

\newcommand{\inc}[1]{\chi^{#1}}
\newcommand{\ktp}{C_3^{+}}
\newcommand{\oktp}{odd-$\ktp$\xspace}
\newcommand{\oktpf}{odd-$\ktp$-free\xspace}
\renewcommand{\phi}{\varphi}
\newcommand{\phib}[1]{\overline{\varphi}(#1)}

\newcommand{\peterm}{\mathbf{T}}

\newcommand{\ci}[1]{\chi'(#1)}
\newcommand{\fci}[1]{\chi'_f(#1)}

\newcommand{\toskf}{totally odd subdivision of $K_4$\xspace}
\usepackage{caption}

\title{Ear-decompositions and the complexity \\ of the matching polytope}
\author{Yohann \sc{Benchetrit}\thanks{Supported by LabEx PERSYVAL-Lab (ANR-11-LABX-0025)}  \and Andr\'as \sc{Seb\H{o}}*}

\date{}
\begin{document}
\maketitle

\begin{abstract}
The complexity of the matching polytope of graphs may be measured with the maximum
length $\beta$ of a starting sequence of odd ears in an ear-decomposition. Indeed, a theorem of
Edmonds and Pulleyblank shows that its facets are defined by 2-connected factor-critical graphs, which have an odd ear-decomposition (according to a theorem of Lov\'asz). 

In particular, $\beta(G)\leq 1$ if and only if the matching polytope of the graph $ G $ is completely described by non-negativity, star and odd-circuit inequalities. This is essentially equivalent to the h-perfection
of the line-graph of $ G $, as observed by Cao and Nemhauser.

The complexity of computing $\beta$ is apparently not known. We show that deciding whether $\beta\leq 1$ can be executed efficiently by looking at any ear-decomposition starting with an odd circuit and performing basic modulo-2 computations. Such a greedy-approach is surprising in view of the complexity of the problem in more special cases by Bruhn and Schaudt, and it is simpler than using the Parity Minor Algorithm.

Our results imply a simple polynomial-time algorithm testing h-perfection in line-graphs (deciding h-perfection is open in general). We also generalize our approach to binary matroids, and show that computing $\beta$ is a Fixed-Parameter-Tractable problem (FPT).

\medskip\noindent
\noindent{{\bf keywords:} ear-decomposition, $2$-connected graphs, odd circuits, cycle space, stable sets, matchings, polytopes, h-perfect graphs, binary matroids, edge-colorings.}

\end{abstract}

\section{Introduction}\label{sec-intro}
In this paper, we only consider finite undirected graphs. They can have multiple edges but no loops. A graph is \emph{simple} if it does not have a pair of parallel edges. We say that a graph $G$ \emph{contains} a graph $H$ if $H$ is a subgraph of $G$.

A \emph{stable set} (resp. \emph{clique}) of a graph is a set of pairwise non-adjacent (resp.  adjacent) vertices. The \emph{chromatic number} of $G$ is the smallest number of stable sets covering $\vts{G}$.

A graph is \emph{perfect} if the chromatic number of each induced subgraph $H$ is equal to the largest cardinality of a clique of $H$. Finding a maximum-weight stable set (or clique) and computing the chromatic number can be carried out in polynomial-time in the class of perfect graphs \cite{Groetschel1988}, whereas these problems are \np-hard in general \cite{Karp1972}.
Besides, deciding whether a graph is perfect can be done efficiently \cite{Cornuejols2003}.

The \emph{incidence vector} of a subset $X$ of a finite set $Y$, denoted $\inc{X}$ is the element of $\set{0,1}^{Y}$ defined for each $y\in Y$ by: $\inc{X}(y)=1$ if and only if $y\in X$.
The \emph{stable set polytope} of a graph $G$ is the convex hull of the incidence vectors of the stable sets of $G$. The \emph{non-negativity inequalities} $x_v\geq 0$ (for each $v\in\vts{G}$) and the \emph{clique inequalities} $\sum_{v\in K}x_v\leq 1$ (for each inclusion-wise maximal clique $K$ of $G$) always define facets of $\stab{G}$ \cite{Padberg1973}.

Results of Lov\'asz \cite{Lovasz1972} and Fulkerson \cite{Fulkerson1972} imply, as stated by Chv\'atal \cite{Chvatal1975}:

\begin{thm}[\cite{Chvatal1975}]\label{thm-LFC}
A graph is perfect if and only if its stable set polytope is described by non-negativity and clique inequalities.
\end{thm}

A \emph{circuit} of $G$ is a 2-regular connected subgraph of $G$, and it is \emph{odd} if it has an odd number of edges. An \emph{odd-circuit inequality} of $G$ is of the form $\sum_{v\in\vts{C}}x_v\leq\frac{\card{\vts{C}}-1}{2}$, where $C$ is an odd circuit of $G$. It is obviously satisfied by every point of $\stab{G}$.

\paragraph{H-perfect graphs } A graph is \emph{h-perfect} if its stable set polytope is completely described by non-negativity, clique and odd-circuit inequalities. It is \emph{t-perfect}\footnote{h is for hole, and t is for its french translation "trou"} if it furthermore does not contain a clique of size 4.

\myref{thm-LFC} easily implies that perfect graphs are h-perfect. Another wide class of t-perfect graphs is obtained by excluding non-t-perfect subdivisions of $K_4$ as subgraphs \cite{Gerards1998} (these subdivisions are characterized in \cite{Barahona1994}).
Besides, Gr\"otschel, Lov\'asz and Schrijver \cite{Grotschel1986} proved that \emph{a maximum-weight stable set of an h-perfect graph can be found in polynomial-time}.

On the other hand, the computational complexity of testing t-perfection is unknown. T-perfection belongs to co-\np \cite{Schrijver2003} but no combinatorial certificate of t-imperfection is available.

The \emph{line graph} of a graph $G$, denoted $L(G)$, is  the simple graph whose vertex-set is $\eds{G}$ and whose edge-set is the set of pairs of incident edges of $G$; $G$ is called a \emph{source graph} of $L(G)$. A graph is \emph{claw-free} if it does not have an induced subgraph isomorphic to the complete bipartite graph $K_{1,3}$. Claw-free graphs form a proper superclass of line-graphs.
Bruhn and Schaudt  proved:

\begin{thm}[Bruhn, Schaudt \cite{Bruhn2014}]\label{thm-Bruhn-Schaudt}
T-perfection can be tested in polynomial-time in the class of claw-free graphs.
\end{thm}

In this paper, we solve this recognition problem for h-perfect line-graphs. These are essentially more general than t-perfect line graphs. Indeed, the maximum degree of source graphs of $t$-perfect line graphs is at most $3$, and their triangles cannot contain parallel edges whereas the source graphs of h-perfect graphs may have arbitrary high degree, and triangles with many parallel edges.

Besides the generalization we found a simple elementary treatment of the subject using ear-decompositions of $2$-connected graphs and related  mod~2 properties of the cycle-space which turns out to be interesting in its own sake. 

The first step was made by Cao and Nemhauser \cite{Cao1998} translating Edmonds and Pulleyblank's \cite{Edmonds1974} complete description of the matching polytope into the line graph.

A \emph{totally odd subdivision} of a graph $H$ is obtained  by replacing each edge $e$ of $H$ with a path having an odd number of edges joining the ends of $e$, such that paths corresponding to distinct edges do not share inner vertices.
Let $ \ktp $ denote the graph obtained from the triangle $K_3$ by adding a single parallel edge. An \emph{\oktp} is a totally odd subdivision of $\ktp$  (they are also called \emph{skewed thetas} \cite{Bruhn2014}).
An \oktp is \emph{strict} if it is not $ \ktp $ itself.

\begin{thm}[Cao, Nemhauser \cite{Cao1998}]\label{thm-CaoNem}
For every graph $H$, the following statements are equivalent:
\begin{itemize}
\item [i)] $L(H)$ is h-perfect,
\item [ii)] $H$ does not contain a strict \oktp.
\end{itemize}
\end{thm}
This extends a previous characterization and algorithm by Trotter \cite{Trotter1977} for perfect line-graphs.
Since \emph{deciding whether a graph $G$ is a line-graph (and building a graph $H$ such that $G=L(H)$ if it exists) can be done in polynomial-time} \cite{Roussopoulos1973}, testing h-perfection in line-graphs reduces to detecting strict \oktp subgraphs.

\begin{figure}[h]
	\centering
	\includegraphics[scale=1.3]{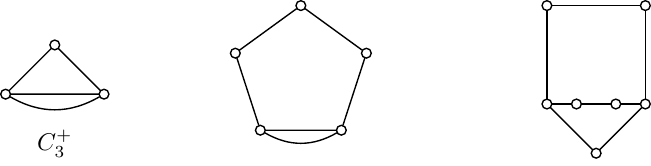}
	\caption{$C_3^{+}$ and two strict \oktp graphs}\label{line-fig-C3p}
\end{figure}

\paragraph{Detecting \oktp subgraphs} Kawarabayashi, Reed and Wollan \cite{Kawarabayashi2011} (and independently Huynh \cite{Huynh2009}) proved the following:

\begin{thm}[\cite{Kawarabayashi2011,Huynh2009}]\label{thm-Kawa}
Let $H$ be a graph. Deciding whether a graph contains a totally odd subdivision of $H$ can be done in polynomial-time.
\end{thm}

A graph is \emph{\oktpf} if it does not contain an \oktp.

Even though \myref{thm-Kawa} can detect an \oktp in an arbitrary, not necessarily simple graph, this is not exactly the obstruction for h-perfection in line graphs according to \myref{thm-CaoNem}; we have to deal only with {\em strict } \oktp.  However, in \myref{sec-algoCfreePlus} we observe that the non-simple strict \oktpf subgraphs can be separately and easily detected. Hence, it only remains detecting an \oktpf in simple graphs, in which all \oktp subgraphs are strict.

Hence \myref{thm-Kawa} will already easily imply:

\begin{restatable}{thm}{HPLine}\label{thm-HPLine}
H-perfection can be tested in polynomial-time in the class of line-graphs.
\end{restatable}

\noindent
We do not know whether h-perfection can be also tested efficiently in the larger class of claw-free graphs.

\myref{thm-Kawa} is built upon elaborated techniques of the Graph Minor Project of Robertson and Seymour and is oriented towards generality. This suggests the search for a more adapted algorithm testing whether a graph is \oktpf.  In this direction, Bruhn and Schaudt \cite{Bruhn2014} provided a direct solution for graphs with maximum degree $3$.

The central contribution of this paper is \emph{a simple  polynomial-time algorithm for the recognition of \oktpf graphs} relying  on a  combinatorial \emph{good characterization theorem for the existence of \oktp in graphs} (that is an \np characterization of \oktpf graphs).   This theorem and its proof are elementary, they avoid  Graph Minors and use the cycle space of a graph instead.

\paragraph{Matroid generalization} A matroid is \emph{binary} if it is the column-matroid of a matrix with coefficients in the field of two elements. The class of binary matroids contains graphic and co-graphic matroids (see \cite{Oxley1992}).

We generalize our approach (algorithms included) to binary matroids. It is surprising that we do not even need ear-decompositions to deal with this more general case, and use only a direct consequence of a theorem of Lehman \cite{Lehman1964}. In particular, this binary generalization provides a different proof and algorithm for the graphic case. Still the graphic case is treated apart, as ear-decompositions show a link with factor-critical subgraphs and h-perfection of line graphs (see also the last paragraph of \myref{sec-matro}).

Complexity of algorithms whose input includes matroids is often measured using the number of required calls to an \emph{independence oracle} (or any other polynomially-equivalent oracle, see \cite{Frank2011}), that is an algorithm testing whether a subset of the ground-set is independent. 

An \emph{\oktp} of a matroid $M$ is a restriction of $M$ which is isomorphic to the circuit  matroid of an \oktp. A matroid is \emph{\oktpf} if it does not have an \oktp. We prove:

\begin{restatable}{thm}{MatroAlgo}\label{thm-MatroAlgo}
Deciding whether a binary matroid $M$ is \oktpf or finding an \oktp of $M$ can be done in polynomial-time using an independence-oracle.
\end{restatable}

\noindent
Our algorithm cannot be directly extended to non-binary matroids and we do not know the complexity of the problem in arbitrary matroids.

\paragraph{Complexity of the matching polytope} We present a new combinatorial parameter motivated by the nice structure of \oktpf graphs and related to the matching polytope.

\begin{restatable}{defin}{DefBeta}\label{def-Beta}
For each graph $G$, let $\beta(G)$ denote the largest  integer $k$ such that $G$ contains a graph $H$ having an odd ear-decomposition with $k$ ears.
\end{restatable}

For example, a graph $G$ is \oktpf if and only if $\beta(G)\leq 1$.
We observe that \myref{thm-Kawa} easily implies: \emph{for each fixed $k$, deciding whether a graph $G$ satisfies $\beta(G)= k$ can be done in polynomial-time}. In other words:

\begin{restatable}{thm}{FPT}\label{thm-FPT}
Determining $\beta$ is a Fixed-Parameter-Tractable problem.
\end{restatable}
We do not even know whether the property $\beta(G)\geq k$ (for each graph $G$ and integer $k$) admits a co-\np-characterization, while the definition clearly shows that it belongs to \np.

The \emph{matching polytope} of a graph is the convex-hull of the incidence vectors of its matchings (a \emph{matching} is a set of pairwise non-incident edges). In other words, it is the stable set polytope of its line graph.
Results of Edmonds, Pulleyblank \cite{Edmonds1974} and Lov\'asz \cite{Lovasz1972a} show that $\beta(G)$ can be used as a parameter to separate on, for questions related to the matching polytope (see the following paragraph on edge-colorings).

The largest number of odd ears in an ear-decomposition of a $2$-connected graph, denoted $\overline{\phi}$, was introduced and studied by Frank in \cite{Frank1993} (in the equivalent form of the smallest number of even ears). We show \emph{a family of graphs for which $\beta=2$  while $\overline{\phi}$ is arbitrarily large}.

\paragraph{$ \beta $ and edge-colorings} The \emph{chromatic index} of a graph $G$, denoted $ \ci{G} $, is the smallest cardinality of a family of  matchings $\mathcal{F}$ such that each edge of $G$ belongs to at least one element of $\mathcal{F}$.
The \emph{fractional chromatic index} of $G$, denoted $\fci{G}$, is the minimum value of $\lambda_1+\cdots+\lambda_k$ with $\lambda_1,\ldots,\lambda_k\in\rsetp$ such that there exist matchings $M_1,\ldots,M_k$ of $G$ satisfying, for each edge $e$ of $G$:
$\sum_{i\in\cro{k}\colon e\in M_i}\lambda_i\geq 1$.

\begin{figure}[h]
	\centering
	 \includegraphics[scale=1]{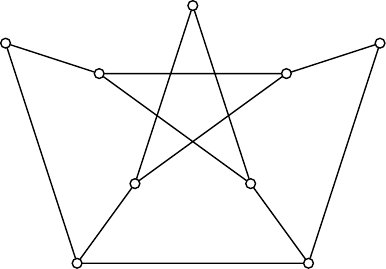}
	\caption{the Petersen graph minus a vertex}\label{line-fig-Pmv}
\end{figure}

It is well-known that the chromatic index of a graph cannot always be  obtained by rounding-up the fractional chromatic index, the smallest known example being the \emph{Petersen graph minus a vertex} (denoted $\peterm$ and shown in \myref{line-fig-Pmv}). Indeed $\ci{\peterm}=4$, whereas Edmonds' description of the matching polytope \cite{Edmonds1965} easily shows that $\fci{\peterm}=3$.

For each graph $G$, let $\hat{G}$ denote the underlying simple graph of $G$. Benchetrit proved:

\begin{thm}[Benchetrit \cite{Benchetrit}]\label{thm-Ben}
Each graph $G$ with $\beta(\hat{G})\leq 1$ satisfies $ \ci{G}=\ceil{\fci{G}} $.
\end{thm}

We conjecture that this result can be extended as follows:

\begin{restatable}{conj}{ConjEcol}\label{conj-Ecol}
Each graph $G$ with $\beta(\hat{G})\leq 3$ satisfies $ \ci{G}=\ceil{\fci{G}} $.
\end{restatable}
The bound 3 would be best possible. Indeed, $\beta(\peterm)=4$ (see \myref{sec-beta}). By results of Baum and Trotter \cite{Baum1981}, this conjecture would imply that the matching polytope $P$ of a graph $G$ with $\beta(\hat{G})\leq 3$ has the \emph{integer decomposition property}: each integral vector of the form $kx$ with $x\in P$ is the sum of $k$ integral vectors of $P$. 

Furthermore, \myref{conj-Ecol} would yield a new case of conjectures of Goldberg \cite{Goldberg1973} and Seymour \cite{Seymour1979} which state:

\begin{conj}[Goldberg \cite{Goldberg1973}, Seymour \cite{Seymour1979a}]\label{line-conj-GoldSeym}
Each graph $G$ satisfies $\chi'(G)\leq\ceil{\chi_f'(G)}+1$.
\end{conj}

Shepherd and Kilakos \cite{Kilakos1996} conjecture that \emph{every graph $G$ which does not have $\peterm$ as a minor satisfies $ \chi'(G)=\ceil{\chi_f'(G)} $}. This would imply that the matching polytope of such graphs has the \idp.
This and \myref{conj-Ecol} do not clearly imply one another. Indeed, it is easy to find graphs without $\peterm$ as a minor and with an arbitrarily large value of $\beta$. Also, the graph obtained from $\peterm$ by subdividing each edge exactly once is bipartite (that is $\beta=0$) and has obviously  $\peterm$ as a minor.

\paragraph{Related works} Cao's thesis \cite{Cao1995} suggests that totally odd subdivisions of $K_4$ are involved in deciding whether a graph is \oktpf. This led us to show \emph{a simple efficient algorithm using $\overline{\phi}$ to detect totally odd subdivisions of $K_4$ in \oktpf graphs}.

The currently known algorithms for detecting such subdivisions in arbitrary graphs are not elementary: \myref{thm-Kawa} directly provides one, and Kawarabayashi, Li and Reed \cite{Kawarabayashi2010} gave a simpler and more adapted algorithm. Both use techniques of the Graph Minor Project. Our simplification for \oktpf graphs is rather specific and does not directly extend to larger values of $\beta$. This does not exclude a possible use of  $\beta$ for  a more general algorithm.

We end the paper with a  review of the results of \cite{Cao1995} concerning \oktpf graphs and observe that some of the statements are incorrect. In particular, the construction procedure given for simple \oktpf graphs does not work.

\textbf{Outline} In \myref{sec-algoCfreePlus}, we first observe that any efficient algorithm deciding whether a simple graph is \oktpf can be used as a black-box to test whether a line-graph is h-perfect in polynomial-time. Hence, we already obtain \myref{thm-HPLine} from \myref{thm-Kawa}.
Then, we prove our characterization of \oktpf graphs in terms of cycle bases and use it to build our  efficient algorithm testing whether a graph (simple or not) is \oktpf.
We extend these ideas to binary matroids in \myref{sec-matro} and prove \myref{thm-MatroAlgo}.

In \myref{sec-beta}, we explain the relation of $\beta$ with the matching polytope and observe that \myref{thm-Kawa} easily implies \myref{thm-FPT}. We also show that $\beta$ and the largest number $\overline{\phi}$ of odd ears in an ear-decomposition need not to be close in general.

We use Frank's algorithm to compute $\overline{\phi}$ \cite{Frank1993} in \myref{sec-SubdivAndCao} to detect totally odd subdivisions of $K_4$ in \oktpf graphs, and finally discuss the related results of \cite{Cao1995} in \myref{subsec-Cao}.

\subsection{Definitions and preliminary results}\label{sec-notations}

For a non-negative integer $k$, we write $\cro{k}$ for the set of integers $1,\ldots,k$.
Let $G$ be a graph and $v$ be a vertex of $G$. The \emph{degree} of $v$ in $G$ is the number of edges incident to $v$ and $\Delta(G)$ is the largest degree of a vertex of $G$.
We write $N_G(v)$ for the set of neighbors of $G$.

For a subset $X$ of $\vts{G}$ or $\eds{G}$, let $G-X$ denote the graph obtained from $G$  by deleting each element of $X$ from $G$. If $X=\set{x}$, we simply write it $G-x$.

A subgraph of $G$ is \emph{induced} if it is obtained from $G$ by deleting vertices. For two graphs $G_1$ and $G_2$, we write $G_1\cup G_2$ for the graph $(\vts{G_1}\cup\vts{G_2},\eds{G_1}\cup\eds{G_2})$.

A \emph{circuit} is a 2-regular connected graph, and a \emph{path} is a circuit minus an edge. So a path has two different vertices of degree one, called its \emph{ends}.
For sets $X,Y\se\vts{G}$, an \emph{$\set{X,Y}$-path} of $G$ is a path joining a vertex of $X$ to a vertex of $Y$. If $X=\set{x}$ and $Y=\set{y}$, then we refer to it as an \emph{$xy$-path} of $G$. Two paths are \emph{inner-disjoint} if they do not share vertices other than their ends.

The \emph{length} of a path (or circuit) is its number of edges. A path (or circuit) is \emph{odd} if it its length is odd, and it is \emph{even} otherwise. A graph is \emph{bipartite} if it does not have an odd circuit.

A connected graph $G$ with at least 3 vertices is \emph{$2$-connected} if $G-v$ is connected for all $v\in V(G)$. A \emph{block} of a graph $G$ is a maximal $2$-connected subgraph or a bridge of $G$ (a \emph{bridge} is the pair of ends of an edge $e$ such that deleting $e$ and all its parallel edges increases the number of components). 

An \emph{ear} of a subgraph $H$ of a graph $G$ is a path of $G$ which has exactly his two different ends in $G$.
An \emph{ear-decomposition} of a graph $G$ is a sequence $(C,P_1,\ldots,P_k)$ of a circuit $C$ and paths $P_i$  such that 
$G=C\cup P_1\cup\cdots \cup P_k$ and, for each $i\in\cro{1,\ldots,k}$: $P_i$ is an ear of $C\cup(P_1\cup\cdots\cup P_{i-1})$. The graphs $C,P_1,\ldots, P_k$ are the \emph{ears} of the decomposition (we omit the usual qualifier ``open'', since we consider only open ear-decompositions).  An ear-decomposition is \emph{odd} if all its ears are odd.


%


\begin{thm}[Whitney \cite{Whitney1932},Robbins \cite{Robbins1939}]\label{thm-whitney}
A graph has an  ear-decomposition if and only if it is 2-connected.
\end{thm}
Besides, we use that \emph{all the ear-decompositions of a $2$-connected graph $G$ have the same number of ears, which is $|\eds{G}|-\card{\vts{G}}+1$}. This follows directly from observing that deleting an edge in each ear of an ear-decomposition of $G$ yields a spanning tree of $G$.

Hence, we may speak of the \emph{number of ears} of a $2$-connected graph (also known as the \emph{cyclomatic number} of the graph).

We frequently use Menger's theorem stating that \emph{for each $2$-connected graph $G$ and each sets $S,T\se\vts{G}$ of cardinality at least 2, there exist two vertex-disjoint $\set{S,T}$-paths}, and that those paths can be found in polynomial-time (see \cite{Schrijver2003}, and \cite{Tholey2006} for recent developments).

\begin{prop}\label{fonda-prop-EarCompletion}
Let $G$ be a $2$-connected graph. Each ear-decomposition of a $2$-connected subgraph  of $G$ can be completed into an ear-decomposition of $G$.
\end{prop}

Several polynomial-time algorithms are available for finding (or completing) an ear-decomposition of a $2$-connected graph (see \cite{Schmidt2013} for a recent example). Also, parallel algorithms were given by Lov\'asz \cite{Lovasz1985} and Miller, Ramachandran \cite{Miller1986}.

Finally, we will frequently use the following easy fact: \emph{if $G$ is a $2$-connected non-bipartite graph, then $G$ contains both odd and even $uv$-paths for each pair of vertices $u$ and $v$ of $G$}. This follows directly by applying Menger's theorem to find two vertex-disjoint paths joining respectively $u$ and $v$ to an arbitrary odd circuit of $G$. This determines the easy and well-known characterization of deciding the existence and finding (in polynomial-time)  {\em a path of given parity between any two vertices of a graph} in terms of its blocks.

\section{A greedy algorithm for recognizing \oktpf graphs}\label{sec-algoCfreePlus}

In this section we prove the main results of the paper. 
We first observe that detecting non-simple \oktp subgraphs can be carried out straightforwardly, and that after filtering these the problem of finding a strict \oktp is equivalent to detect an \oktp in the underlying simple graph. Hence, detecting a (strict or not) \oktp of a graph $H$ turn out to be the only essential difficulty in testing the h-perfection of $L(H)$.
Then we show that the existence of an \oktp subgraph can be readily decided from {\em any arbitrary ear-decomposition of $G$ starting with an odd circuit} in a simple, elementary way (based on a good characterization of \oktpf graphs, see \myref{subsec-circSpace}). This is a surprising result in view of complicated previous partial solutions or the use of the graph minor theory of Robertson and Seymour.
Finally, in \myref{sec-matro} we generalize the results to binary matroids.

\medskip




Since recognizing line-graphs (and building a corresponding source graph if it exists) can be done efficiently \cite{Roussopoulos1973}, \myref{thm-CaoNem} shows that deciding h-perfection in line-graphs reduces to detecting  strict \oktp subgraph in the source graph.

Clearly, {\em the only strict \oktp graphs which are not simple consist of an odd circuit of length at least $5$ with two neighboring vertices $u,v$ joined by two parallel edges.}   Such an odd circuit can be easily detected in polynomial time :

\begin{prop}\label{line-prop-strict}
Let $G$ be a graph, $u$ and $v$ be two vertices of $G$. Finding an even $uv$-path of length at least $4$ or certifying its non-existence can be done in polynomial time.
\end{prop}

\begin{dem}
There exists an even $uv$-path of length at least $4$ if and only if there exist $a\in N_G(u)$ and an odd $av$-path in $(G - u) - av$ (we mean that all edges whose ends are $a$ and $v$ are deleted).  We then use  that odd paths between two vertices of a graph can be found or proved not to exist in polynomial time (see last paragraph of Section~\ref{sec-notations}).
\end{dem}

Note that {\em \myref{line-prop-strict} contains the problem of detecting an odd circuit of length at least $5$ through a given edge $uv$}, that is an odd hole containing a given vertex in the line graph. An efficient algorithm for this problem for the considerably larger class of claw-free graphs is given in \cite{VantHof2012}. This is an \np-complete problem in graphs in general \cite{Bienstock1991}.

Using \myref{line-prop-strict} for  all $u$ and $v$ with at least two parallel edges between them means detecting non-simple strict \oktp graphs or certifying that they do not exist. 
It remains to detect simple \oktp subgraphs or proving that the input graph is \oktpf, which is a priori more difficult (see \myref{thm-CaoNem} and \myref{thm-Kawa}). We solve this task in a simple self-contained way in \myref{subsec-circSpace}.  



Let us note that in the particular case of graphs of maximum degree $3$, Bruhn and Schaudt also provided an algorithm detecting \oktp which is elementary and avoids Graph Minors. 

\subsection{A binary characterization of \oktpf graphs}\label{subsec-circSpace}
We write $\ftwo$ for the field of two elements.
Let $G$ be a graph. Clearly, the sum in the vector space $\ftwo^{\eds{G}}$ of the incidence vectors of $F_1\se\eds{G}$ and $F_2\se\eds{G}$ is the incidence vector of the symmetric difference $F_1\Delta F_2$.

A \emph{cycle} is the union of edge-disjoint circuits of $G$ (identified to their edge-sets); equivalently, it is a subgraph with all degrees even.   The \emph{cycle space of $G$}, denoted $\cspace{G}$, is the subspace of the vector space $\ftwo^{\eds{G}}$ consisting of the incidence vectors of cycles. It is spanned (over $\ftwo$) by  the incidence vectors of the circuits of $G$. The rank of $\cspace{G}$ is $|E(G)| - |V(G)| +1$ if $G$ is connected and is the {\em cyclomatic number} of $G$. A well-known class of bases of $\cspace{G}$ is obtained as follows: take any fixed spanning tree $T$ of $G$ and for each $e\in\eds{G}\sm\eds{T}$, let $C_e$ be the unique circuit of $T+e$. It is straightforward ot check that the incidence vectors of circuits $C_e$ obtained form a basis of $\cspace{G}$. 

A \emph{cycle basis} of $G$ is a set of cycles whose incidence vectors form a basis of $\cspace{G}$ (over $\ftwo$). If all members of a cycle basis are circuits, then we call it a \emph{circuit basis}. 

A cycle of a graph is \emph{odd} if it has an odd number of edges, and a cycle basis of a graph is \emph{odd} if all its elements are odd. An odd cycle basis of a graph is \emph{totally odd} if its odd cycles pairwise-intersect in an odd number of edges. 
For example, each set of 3 circuits of a totally odd subdivision  of $K_4$ form a totally odd circuit basis of $K_4$.

In this section, we prove the following characterization of \oktpf graphs and use it to build our algorithm for the recognition of these graphs.

Since an \oktp is 2-connected and non-bipartite, we need only to consider 2-connected non-bipartite graphs.

\begin{restatable}{thm}{CaracTOddBasis}\label{line-thm-CaracTOddBasis}
Let $G$ be a 2-connected non-bipartite graph. The following statements are equivalent:
\begin{itemize}
\item [(i)] $G$ is \oktpf,
\item [(ii)] $G$ has a totally odd circuit basis,
\item [(iii)] each odd cycle basis of $G$ is totally odd.
\end{itemize}
\end{restatable}

We first state a few results needed for proving this theorem.

Cao's thesis \cite{Cao1995} shows that \emph{the odd circuits of a 2-connected \oktpf simple graph pairwise-intersect in an odd number of edges}.
We first observe that this property characterizes 2-connected \oktpf graphs:

\begin{lem}\label{line-prop-oddInter}

Let $G$ be a 2-connected graph. Then  $G$ contains a \oktp if and only if it has two odd circuits meeting in an even number of elements.

Furthermore, from two such odd circuits of a 2-connected graph an \oktp can be constructed in polynomial-time
\end{lem}

\begin{dem}
Clearly, an \oktp has exactly two odd circuits which have an even number of common edges and thus (ii)=>i).

Conversely, suppose that $G$ has odd circuits $C_1$ and $C_2$ such that $\card{\eds{C_1}\cap\eds{C_2}}$ is even. We show that $G$ contains an \oktp.

First, let us assume that $\card{\vts{C_1}\cap\vts{C_2}}\leq 1$.
Since $G$ is $2$-connected, Menger's theorem shows that there exist two vertex-disjoint $ \set{\vts{C_1},\vts{C_2}} $-paths $P$ and $Q$ (one may be reduced to a single vertex if $C_1$ and $C_2$ meet).
Let $p$ and $q$ be the respective ends of $P$ and $Q$ on $C_1$ and let $R$ be the unique $pq$-path of $C_1$ whose parity is distinct from $\card{\eds{P}}+\card{\eds{Q}}$. Clearly, $R\cup P\cup Q\cup C_2$ is an \oktp subgraph of $G$.

Now, suppose that $C_1$ and $C_2$ have at least two vertices in common.
Since both circuits are odd and $C_1\neq C_2$, the set $\vts{C_1}\cap \vts{C_2}$ defines a partition of $C_1$ into edge-disjoint paths $P_1,\ldots, P_k$ $(k\geq 1)$ which have exactly their ends in $\vts{C_2}$. Since $\card{\eds{C_1}\cap\eds{C_2}}$ is even and as $C_1$ is odd, at least one of these paths must be odd, say $P_1$, and $C_2\cup P_1$ is an \oktp of $G$.

\end{dem}


The proof is clearly algorithmic.




The following lemma plays a key-role in the proof of \myref{line-thm-CaracTOddBasis}; the fact that we have only circuits in the basis is important ! 

\begin{restatable}{lem}{oddBasisExistence}\label{line-prop-oddBasisExistence}
Each 2-connected non-bipartite graph has an odd circuit basis. Furthermore, such a circuit basis can be found in polynomial-time.
\end{restatable}

\begin{dem}
Let $G$ be a $2$-connected non-bipartite graph and $C$ be an odd circuit of $G$.
By \myref{fonda-prop-EarCompletion}, $G$ has an ear-decomposition $(P_0,P_1,\ldots,P_k)$. Recall from \myref{sec-notations} that $k=\card{\eds{G}}-\card{\vts{G}}+1$, which is the cyclomatic number of $G$.

For each $i\in\cro{k}$, the graph $C\cup P_1\cdots\cup P_{i-1}$ is $2$-connected and non-bipartite, so it contains a path $Q_i$ which joins the ends of $P_i$ and such that the circuit $P_i\cup Q_i$ is odd. It is straightforward to check that the incidence vectors of the circuits $P_1\cup Q_1,\ldots,P_k\cup Q_k$ are linearly independent. Hence, they form a circuit basis of $G$, which is odd.

An ear-decomposition and the paths $Q_i$ can be computed in polynomial-time (see \myref{sec-notations}).
\end{dem}

\noindent
In general, a cycle basis does not need to contain only circuits: the fact that here it consists only of particular circuits is a key-point of our proof. A relevant property of totally odd bases can be extended to all odd circuits of the graph:

\begin{restatable}{lem}{ToddImpliesOddInter}\label{line-prop-ToddImpliesOddInter}
If a 2-connected graph has a totally odd cycle basis, then any odd cycles of the graph intersect on an odd number of common edges.
\end{restatable}

\begin{dem}
Let $\cdot$ denote the standard bilinear form on $\ftwo^{\eds{G}}$. That is, for subsets $F_1$ and $F_2$ of $\eds{G}$: $\chi^{F_1}\cdot \chi^{F_2}$ is equal to 1 if $\card{F_1\cap F_2}$ is odd, and 0 otherwise. \emph{Until the end of this proof, all equalities take place in $\ftwo^{\eds{G}}$.}

Suppose that $G$ has a totally odd cycle basis $\mathcal{B}$ and let $C_1$ and $C_2$ be odd cycles of $G$. We show that $\chi^{\eds{C_1}}\cdot\chi^{\eds{C_2}}=1$, as stated.

Since $\mathcal{B}$ is a cycle basis of $G$, there exists $\mathcal{B}_1\se\mathcal{B}$ and $\mathcal{B}_2\se\mathcal{B}$ such that:
\begin{center}
$\chi^{\eds{C_1}}=\displaystyle\sum_{C\in \mathcal{B}_1}\chi^{\eds{C}}$ and $\chi^{\eds{C_2}}=\displaystyle\sum_{D\in\mathcal{B}_2}\chi^{\eds{D}}$.
\end{center}

Since $C_1$ and $\mathcal{B}$ are odd, multiplying by the all-1 vector $\one$ on both sides of the first equality yields: $\card{\mathcal{B}}=1$ (that is, $\mathcal{B}_1$ has odd cardinality).
Similarly, $\card{\mathcal{B}_2}=1$. Since $\mathcal{B}$ is totally odd, we obtain by linearity:

\[ \inc{\eds{C_1}}\cdot \inc{\eds{C_2}}=\sum_{C\in\mathcal{B}_1,\,D\in\mathcal{B}_2}\inc{\eds{C}}\cdot\inc{\eds{D}}
=\sum_{C\in\mathcal{B}_1,\,D\in\mathcal{B}_2}1=\card{\mathcal{B}_1}\card{\mathcal{B}_2}=1, \]
and this ends the proof of the proposition.
\end{dem}

We now prove \myref{line-thm-CaracTOddBasis} using those preliminary results:

\begin{dem}[of \myref{line-thm-CaracTOddBasis}]
We first show that i)=>ii).
Suppose that $G$ is \oktpf. Since $G$ is 2-connected and non-bipartite, \myref{line-prop-oddBasisExistence} shows that $G$ has an odd cycle basis $\set{C_1,\ldots,C_k}$ such that each $C_i$  is a circuit ($i\in\cro{k}$).

As $G$ is \oktpf, \myref{line-prop-oddInter} shows that the odd circuits $C_1,\ldots, C_k$ pairwise-intersect in an odd number of edges. Therefore, the basis $\set{C_1,\ldots,C_k}$ is totally odd.

The implication ii)=>iii) straightforwardly follows from \myref{line-prop-ToddImpliesOddInter}. We now show iii)=>i).

Suppose that each odd cycle basis of $G$ is totally odd.
Since $G$ is 2-connected and non-bipartite, \myref{line-prop-oddBasisExistence} shows that $G$ has an odd cycle basis $\mathcal{B}$. By assumption, $\mathcal{B}$ is totally odd. Hence, \myref{line-prop-ToddImpliesOddInter} implies that odd cycles, and in particular odd circuits, pairwise-intersect in an odd number of edges. By \myref{line-prop-oddInter}, this shows that $G$ is \oktpf.
\end{dem}

Clearly, this proof of \myref{line-thm-CaracTOddBasis} provides an algorithm deciding whether a graph is \oktpf: we first build efficiently an odd circuit basis \myref{line-prop-oddBasisExistence}. If there are two odd circuits of the basis having an even number of common edges, we build an \oktp from them. Otherwise, either $G$ is bipartite or any pair of odd circuits in the basis meet in an odd number of elements, certifying that the basis is totally odd and that $G$ is \oktpf (\myref{line-prop-ToddImpliesOddInter}).


\subsection{Extension to binary matroids}\label{sec-matro}
In this section we show that the results of the previous section can be generalized to binary matroids. 
Standard terminology and basic facts related to matroids can be found for instance in \cite{Oxley1992} (binary matroids are treated in Chapter 9 of this book). We consider \emph{loopless} matroids only.

A  matroid is \emph{binary} if it is representable in a linear space over $\ftwo$. It is well-known that \emph{a matroid  is binary if and only if the symmetric difference of any set of circuits of is the union of disjoint circuits}. 

We say that a matroid is an \emph{\oktp} if it is the circuit matroid of an \oktp graph (recall that an \oktp graph is a totally odd subdivision of the graph $C_3^{+}$). An \emph{\oktp of a matroid} $M$ is a restriction of $M$ which is isomorphic to an \oktp, and a matroid is \emph{\oktpf} if it does not have an \oktp.

Connectivity assumptions were important in our treatment for \oktp graphs. It is the same for the proof of our matroid generalization and we thus recall the corresponding notions here. 
Let $M$ be a matroid with ground set $S$. Consider the relation on $S$ defined by: $e,f\in S$ are related if and only if $e=f$ or there exists a circuit containing both $e$ and $f$. It is well-known that this is an equivalence relation, whose classes are called the \emph{blocks} of $M$. A matroid is \emph{connected} if it has at least two elements and only one block. Note that the connectedness of the circuit matroid $M$ of a graph $G$ with at least 3 vertices means the 2-connectedness of $G$ and that the blocks of $M$ correspond to the edge-sets of the blocks of $G$.


The following straightforward characterization of \oktp matroids will be useful:

\begin{prop}\label{prop-equivOddC3p}
Let $M$ be a matroid. The following statements are equivalent: 
\begin{itemize}
\item [i)] $M$ is an \oktp,
\item [ii)] $M$ is the union of two circuits $C_1$ and $C_2$ such that $C_1$ is odd, $\card{C_2\sm C_1}$ is odd and $M$ has exactly three circuits which are: $C_1$, $C_2$ and $C_1\Delta C_2$.
\end{itemize}
\end{prop}

Binary matroids generalize both graphic and co-graphic matroids. We extend the cycle-space approach of \myref{subsec-circSpace} for \oktpf graphs to show an efficient algorithm which tests whether a matroid is \oktpf or finds an \oktp otherwise. The input matroid can be given by a linear representation, but we need only an independence oracle (which is in fact equivalent in terms of algorithmic complexity).

\MatroAlgo


A \emph{cycle} of a matroid is a union of disjoint circuits, and it is \emph{odd} if it has an odd number of elements. It is well-known that, as for graphs, the set of (incidence vectors of) cycles of a binary matroid $M$ with ground set $S$ is a subspace of $\ftwo^{S}$ (this actually characterizes binary matroids \cite[chap. 9 9]{Oxley1992}). It is called the \emph{cycle space} of $M$ and is denoted $\cspace{M}$. Clearly, $\cspace{M}$ is the subspace of $\ftwo^{S}$ spanned by the circuits of $M$ and it is easy to check that \emph{the rank of $\cspace{M}$ is $|S|-r$}, where $r$ is the rank of $M$.

Cycle and circuit bases, odd and totally odd cycle bases of $M$ are defined in the exact same way as for graphs (see \myref{sec-algoCfreePlus}).

The following lemma generalizes \myref{line-prop-oddBasisExistence}. The main technical difficulty is to show an \oktp from two given disjoint odd circuits, without the availability of Menger's theorem in graphs.

\begin{lem}\label{matrline-prop-oddInter}
	A connected matroid has an \oktp  if and only if it has two odd circuits which meet in an even number of elements.
		
	
	Furthermore, from two such odd circuits  an \oktp can be constructed in polynomial-time.
\end{lem}

\begin{dem}
Clearly, an \oktp has exactly two odd circuits which have an even number of common elements. 

Conversely, we first show the following:

\medskip\noindent
\textsc{Claim.} \emph{If $C_1$ and $C_2$ are two circuits of a matroid such that $C_1\cap C_2\neq \vn$, $C_1$ is odd and $\card{C_1\sm C_2}$ is odd, then: $C_1\cup C_2$ contains an \oktp.}

Indeed, we prove that if $C_1\cup C_2$ is inclusion-wise minimal among all possible choices respecting the assumptions, then $C_1\cup C_2$ is an \oktp.

Let $C\se C_1\cup C_2$ be a circuit which is neither $C_1$ nor $C_2$ (such a circuit must exists since $C_1$ and $C_2$ meet). We will show that $C=C_1\Delta C_2$, and this and \myref{prop-equivOddC3p} will imply the claim.
Clearly, $C$ must meet both $C_1$ and $C_2$. 

If $\card{C\sm C_2}$ is even, then since $M$ is binary the set $C\Delta C_2$ is the union of disjoint circuits. Since $\card{C_2\sm C_1}$ is odd, one of them, say $C'$, is such that $\card{C'\sm C_1}$ is odd. Hence the pair $(C_1,C')$ satisfies the assumptions of the claim and minimality shows $C'\sm C_1=C_2\sm C_1$, that is $C'\Delta C_2\se C_1$. Since $M$ is binary, this implies that $C'=C_1\Delta C_2=C\Delta C_2$ and thus $C=C_1$: a contradiction.

Therefore, we may assume that $\card{C\sm C_2}$ is odd. Minimality then shows that $C\sm C_1=C_2\sm C_1$. Since $M$ is binary, this implies $C\Delta C_2= C_1$ and we are done.

\medskip\noindent

We now use the claim to prove the lemma. Let $C_1$ and $C_2$ be two odd circuits of $M$ meeting on an even number of elements. Clearly, the claim yields an \oktp if $C_1\cap C_2\neq \vn$ so we may assume the contrary.

Since $M$ is connected, it has a circuit meeting both $C_1$ and $C_2$ and we may consider such a circuit $C$ with $C\sm C_2$ inclusion-wise minimal. 

The set $C\Delta C_2$ is a circuit: indeed since $M$ is binary, $C_1\Delta C_2$ must contain a circuit $C'$ which meets $C_1$ and $C_2$. The minimality of $C\sm C_2$  shows that $C'=C\Delta C_2$, and therefore $C\Delta C_2=C'$ as required.

Both $C$ and $C\Delta C_2$ meet $C_1$ and, since $C_2\sm C_1$ is odd, one of them has an odd number of elements outside of $C_1$. Therefore we may apply the claim again to obtain an \oktp of $M$.

\end{dem}

\medskip
We now prove a generalization of \myref{line-prop-oddBasisExistence} to binary matroids, that makes possible to extend all the results. Surprisingly, we do not need the generalization of ear-decompositions to matroids \cite{Coullard1996} (see also the last paragraph of this section) to prove this but use only the following straightforward consequence of a result of Lehman instead (\cite[chap. 9.3, exercice 9]{Oxley1992}):

\begin{prop}\label{prop-Lehman}
Each element of a binary connected matroid $M$ belongs to an odd circuit of $M$.
\end{prop}

\begin{lem}\label{lem-oddBasisMatro}
For any connected non-bipartite binary matroid there exists an odd circuit basis  that can be constructed in polynomial time.
\end{lem}

\begin{dem} Let $M$ be a connected binary matroid. Let $M_p$ be the binary matroid obtained by adding successively an all-0 column and an all-1 line to a matrix representation of $M$, and let $p$ be the new element of $M_p$.

Using \myref{prop-Lehman}, it is straightforward to check that $M_p$ is a connected matroid. This implies that we can build greedily a set of circuits $C_1,\ldots, C_k$ of $M_p$ which all contain $p$ and such that for each $i\in\set{1,\ldots,k-1}$: $C_{i+1}\sm C_i\neq\vn$ . It is now straightforward to check that $\set{C_1-p,\ldots,C_k-p}$ is an odd circuit basis of $M$ (and $k$ is the number of elements of $M$ minus its rank).


\end{dem}

Now we can immediately extend \myref{line-prop-ToddImpliesOddInter} and \myref{line-thm-CaracTOddBasis} to binary matroids.  

\begin{restatable}{lem}{MatrToddImpliesOddInter}\label{line-prop-ToddImpliesOddInter}
	If a connected matroid has a totally odd cycle basis, then all its odd cycles pairwise-intersect in an odd number of edges.
\end{restatable}

\begin{thm}\label{thm-EquivMatro}
Let $M$ be a connected non-bipartite binary matroid. The following statements are equivalent:
\begin{itemize}
\item [i)] $M$ is \oktpf,
\item [ii)] $M$ has a totally odd circuit basis,
\item [iii)] each odd cycle basis of $M$ is totally odd.
\end{itemize}
\end{thm}





Since finding the blocks of $M$ can be easily done efficiently, turning the proof of \myref{thm-EquivMatro} into a polynomial-time algorithm testing whether a matroid is \oktpf only requires one more subroutine: deciding efficiently whether a connected binary matroid is bipartite. This can be carried out using the following simple proposition, which generalizes the bipartiteness test of graphs:

\begin{prop}\label{thm-bipartite}
Let $M$ be a connected binary matroid. The following statements are equivalent:
\begin{itemize}
\item [i)] $M$ is bipartite,
\item [ii)] There exists a circuit basis of $M$ containing only even cycles,
\item [iii)] Each cycle basis of $M$ contains only even cycles.
\end{itemize}
\end{prop}

The statements  (i) implies (ii) and (iii) implies (i) are obvious. For the (ii) implies (iii) part, note that the parity of the symmetric difference of two cycles is the mod~$2$ sum of the two parities. This proves the proposition.

It follows that any circuit basis is a good certificate for bipartiteness (so is a well-known fourth equivalent statement as well: {\em the ground set of $M$ is the disjoint union of cocycles}). It also follows that {\em bipartiteness of matroids can be tested in polynomial time.}

We conclude that {\em testing for an \oktp or certifying that the matroid is \oktpf
	 can be solved in polynomial time for binary matroids as well}. 
	 
	 Specialized to graphic matroids, this provides another algorithm testing whether a graph is \oktpf. However, contrarily to the use of ear-decompositions, this alternative approach to building an odd circuit basis is not natural for graphs (as the class of graphic matroids is not closed under the operation $M_p$ used in the proof of \myref{lem-oddBasisMatro}) and it does not directly suggest the relation with the matching polytope discussed in \myref{sec-beta}.

\section{Odd ears and the matching polytope}\label{sec-beta}

In this section, we introduce a new combinatorial parameter, denoted$\beta$, measuring the complexity of facets of the matching polytope and which generalizes \oktpf graphs. We observe that computing it is a Fixed-Parameter-Tractable problem. See \myref{sec-intro} for a useful application of $\beta$ to edge-colorings.

We then discuss the connection of  $\beta$ with the largest number of odd ears in an ear-decomposition (\myref {subsec-Frank})

\subsection{A measure of the complexity of the matching polytope}\label{subsec-matchingPolytope}

We write $\match{G}$ for the \emph{matching polytope} of a graph $G$, that is the convex hull of the incidence vectors of its matchings. For each $v\in\vts{G}$, let $\delta_G(v)$ denote the set of edges incident to $v$.

A graph $G$ is \emph{factor-critical} if for each $v\in\vts{G}$, the graph $G-v$ has a perfect matching. Edmonds and Pulleyblank characterized the facets of the matching polytope. Their results imply:

\begin{restatable}[Edmonds, Pulleyblank \cite{Edmonds1974}]{thm}{EdPull}\label{thm-EdPull}
For every graph $G$:

\[
\match{G}:=\left\{x\in\rset^{\eds{G}}\colon\begin{array}{cc}
  x\geq 0,  &  \\
  \displaystyle \sum_{e\in\delta_G(v)}x_e\leq 1  & \forall v\in\vts{G},  \\
  \displaystyle \sum_{e\in\eds{H}}x_e\leq \frac{|\vts{H}|-1}{2}  & \text{$\forall H$ 2-connected induced} \\
   & \text{factor-critical subgraph of $G$.}
\end{array}\right\}. \]

\end{restatable}

Lov\'asz proved:
\begin{thm}[Lov\'asz \cite{Lovasz1972a,Lovasz2009}]\label{thm-LovaszFCG}
A 2-connected graph is factor-critical if and only if it has an odd ear-decomposition.
\end{thm}

These two results are the main tools for proving \myref{thm-CaoNem} in \cite{Cao1998}.

Together with our results on \oktpf graphs, they motivate us to introduce the following parameter (see \myref{sec-notations} for the definition of the number of ears of a 2-connected graph):

For each 2-connected graph $G$, let $\beta(G)$ denote the maximum number of odd ears starting an ear-decomposition of $G$. 
By \myref{thm-LovaszFCG}, $\beta(G)$ is the largest number of ears of a 2-connected factor-critical subgraph of $G$ and hence this definition of $\beta$ in terms of ears and \myref{def-Beta}  are  equivalent.
Furthermore, \myref{thm-EdPull} shows that $\beta$ can be used as a parameter to separate on, for questions related to the matching polytope (see the paragraph on edge-colorings in \myref{sec-intro}).

Clearly, an \oktp is a 2-connected graph having an ear-decomposition with exactly two ears which are both odd.
Therefore, \emph{a graph $G$ is \oktpf if and only if $\beta(G)\leq 1$}. Besides, \myref{thm-LovaszFCG} states that \emph{a 2-connected graph $G$ is factor-critical if and only if $\beta(G)=\card{\eds{G}}-\card{\vts{G}}+1$}.

The property $\beta\geq k$ obviously belongs to \np. We do not know whether it admits a co-\np characterization.

\begin{quest}
Can $\beta$ be determined in polynomial-time ?
\end{quest}

Let $k$ be a positive integer and $G$ a graph. Clearly, a 2-connected factor-critical graph with $k$ ears is a totally odd subdivision of a graph with at most $ 2k-2 $ vertices of degree at least 3 and at most $ 3k $ edges. Hence checking whether $\beta(G)\geq k$ can be done by enumerating all factor-critical graphs $H$ with $\card{\vts{H}}\leq 2k-2$ and $\card{\eds{H}}\leq 3k$ and use \myref{thm-Kawa} to test whether $G$ contains a totally odd subdivision of $H$. This shows a polynomial-time algorithm deciding $\beta(G)\geq k$ for $k$ fixed, that is:


\FPT*

%
%

We showed in \myref{sec-algoCfreePlus} a simpler efficient algorithm recognizing \oktpf graphs, that is deciding $\beta\leq 1$. We do not know the solution for larger values of $\beta$.

\subsection{Frank's parameter $\phi$. Relation with $\beta$}\label{subsec-Frank}

For a $2$-connected graph $G$, let $\phi(G)$ denote the smallest number of even ears in an ear-decomposition of $G$.
This was introduced by Frank \cite{Frank1993} (for non-necessarily open ear-decompositions), and results of \cite{Frank1993} imply that \emph{an ear-decomposition of a 2-connected graph $G$ with $\phi(G)$ ears can be found efficiently} (see \cite[Section 3]{Cheriyan2001} for a proof).


Let $G$ be $2$-connected, and put:
\[ \phib{G}:=\card{\eds{G}}-\card{\vts{G}}+1-\phi(G).  \]
Since the ear-decompositions of $G$ all have the same number $\card{\eds{G}}-\card{\vts{G}}+1$ of ears (see \myref{sec-notations}), \emph{$\phib{G}$ is the largest number of odd ears in an  ear-decomposition of $G$}.

Each 2-connected graph $G$ obviously satisfies $\phib{G}\geq \beta(G)$. In this section, we show a family of graphs with $\beta=2$ and $\overline{\varphi}$ arbitrarily large.


Let $k\geq 3$ be an integer and $T_1,\ldots,T_k$ be $k$ vertex-disjoint copies of the simple graph obtained from the circuit of length 5 by adding a single edge. Let $v_i$ be the unique vertex of degree 2 in the triangle of $T_i$ and let $u_i$ be one of its neighbors.

\begin{figure}[h]
	\centering
	\includegraphics[scale=1.2]{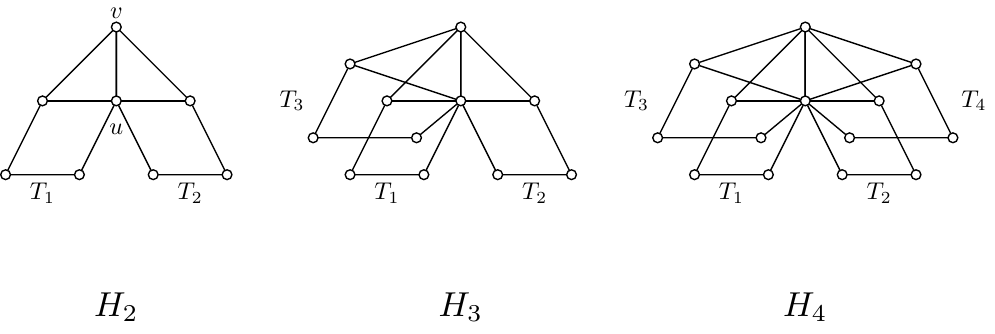}
	\caption{the graphs $H_2$, $H_3$ and $H_4$}\label{line-fig-Hkgraphs}
\end{figure}

Now, let $H_k$ be the graph obtained by identifying all the $v_i$ to a single vertex $v$, all the $u_i$ to a single vertex $u$ and keeping only one copy of the edge $uv$ (see \myref{line-fig-Hkgraphs}).

It is straightforward to check the following:

\begin{prop}\label{line-prop-ArbitrarilyLargediff}
For each $k\geq 2$:
\begin{center}
$\beta(H_k)=2$ and $\phib{H_k}\geq k$.
\end{center}
\end{prop}


In \cite{Frank1993}, Frank showed a min-max theorem for $\phi$ in terms of maximum-cardinality joins: a \emph{join} of a graph $G$ is a set $F\se\eds{G}$ such that each circuit $C$ of $G$ satisfies $|\eds{C}\cap F|\leq|\eds{C}\sm F|$.
We do not know whether a similar min-max result holds for $\beta$.

Even though the much simpler greedy ear-construction of \myref{sec-algoCfreePlus} finally  provided the appropriate answer, the parameter $\phi$ provided a first tool for deciding $\beta\leq 1$ or $\beta\geq 2$ in very particular cases. We sketch in \myref{sec-Appendix} some possibly useful relations.



\appendix
\section{Appendix: subdivisions of $K_4$ and \oktp graphs}\label{sec-Appendix}

\myref{line-prop-ArbitrarilyLargediff} shows that  $\overline{\phi}$ is not really closely related to $\beta$. However, an investigation of their equality may provide new insights. In \myref{sec-SubdivAndCao}, we use Frank's algorithm to compute $\bar \phi$ as a black box to show a rather simple efficient algorithm finding totally odd subdivisions of $K_4$ in \oktpf graphs (and thus a relation between two relevant families of subgraphs). Even though this is irrelevant for the actual discussion on the recognition of \oktpf graphs, it has been a motivation for our work. 
It is the same for related results of Cao's thesis \cite{Cao1995}, from which we got our first inspirations for characterizing h-perfect line-graphs, and we discuss those in \myref{subsec-Cao}.

\subsection{Finding a totally odd subdivisions of $K_4$ in an \oktpf graphs}\label{sec-SubdivAndCao}


Finding a \toskf subgraph is not elementary in general: the simplest algorithm available for their detection in arbitrary graphs uses general techniques of the Graph Minor Project \cite{Kawarabayashi2010}.

Our algorithm is based on the following characterization.
Clearly, we need only to consider simple 2-connected graphs (the following statement is actually false for non-simple graphs in general, as shows the graph obtained by adding two parallel edges to $C_4$).

\begin{restatable}{thm}{KFourfromPhib}\label{line-thm-PhiTwoImpliesK4}
Let $G$ be a 2-connected \oktpf simple graph. The following statements are equivalent:
\begin{itemize}
\item [i)] $G$ does not contain a \toskf,
\item [ii)] $\phib{G}\leq 1$.
\end{itemize}
\end{restatable}

We say that an ear-decomposition of a 2-connected graph $G$ is \emph{optimal} if it has $\phib{G}$ odd ears. As mentioned in \myref{subsec-Frank},
results of \cite{Frank1993} show that 
an optimal ear-decomposition of $G$ can be found in polynomial-time.

Therefore, \myref{line-thm-PhiTwoImpliesK4} directly implies that testing whether an \oktpf graph contains a \toskf can be carried out in polynomial-time. Finding efficiently such a subdivision (if it exists) easily follows from our proof, which is constructive.

An \emph{odd theta} is a graph formed by three inner-disjoint odd paths with the same ends (each path may be reduced to a single edge, see \myref{line-fig-oddThetas}). The first ingredients are the following statements:

\begin{figure}[h]
	\centering
	\includegraphics[scale=1.4]{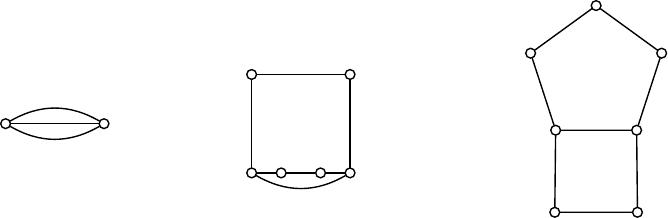}
	\caption{examples of odd thetas}\label{line-fig-oddThetas}
\end{figure}

\begin{restatable}{prop}{oddtheta}\label{line-prop-oddtheta}
Let $G$ be a 2-connected bipartite graph and $H$ a 2-connected subgraph of $G$. If $H$ has an odd  ear in $G$, then each vertex of $G$ belongs to an odd theta subgraph of $G$.
\end{restatable}


\begin{dem}
	Let $P$ be an odd ear of $H$ in $G$.
	We first show that $G$ contains an odd theta. Let $P$ be an odd ear of $H$ and let $u_1$ and $u_2$ be the ends of $P$. 
	
	
	Since $H$ has an ear-decomposition, it is 2-connected. In particular, Menger's theorem shows that $H$ contains two internally vertex-disjoint $u_1u_2$-paths $Q$ and $R$. Since $P$ is odd and $G$ is bipartite, both $Q$ and $R$ are odd. Clearly, $\vts{Q}\cap\vts{R}$ defines a partition of the edge-set of $Q$ into paths $Q_1,\ldots,Q_l$.
	Since $Q$ is odd, one of those paths, say $Q_1$, must be odd. It is easy to check that $R\cup Q_1\cup P$ is an odd theta of $G$.

	Finally, we prove that every vertex of $G$ belongs to an odd theta. Let $T$ be an odd theta of $G$ and let $s\in\vts{G}\sm\vts{T}$. Since $G$ is 2-connected, Menger's theorem  shows that there are two $\set{s,\vts{T}}$-paths $Q_1$ and $Q_2$ whose only common vertex is $s$. A straightforward and short case-checking shows that $Q_1\cup Q_2\cup T$ always has an odd theta containing $s$.
\end{dem}

It is straightforward to convert this proof into a polynomial-time algorithm which finds an odd theta containing a prescribed vertex under the assumptions.

For finishing the proof of \myref{line-thm-PhiTwoImpliesK4} we also need the following:

\begin{restatable}{prop}{findinganoddKfour}\label{line-prop-FindingAtosKfour}
Let $G$ be a 2-connected non-bipartite graph, $C$ an odd circuit of $G$ and $v\in\vts{G}\sm\vts{C}$.
If $G$ contains three inner-disjoint odd  $\set{v,\vts{C}}$-paths, then $G$ contains an \oktp or a \toskf.
\end{restatable}

\begin{dem}
	Let $P_1, P_2, P_3$ be three inner-disjoint odd  $ \set{v,\vts{C}} $-path and let $k:=\card{(\cup_{i=1}^{3}\vts{P_i})\cap\vts{C}}$.
	
	\begin{case} $k=1$. Let $u$ be the unique vertex of $ (\cup_{i=1}^{3}\vts{P_i})\cap\vts{C} $.
		
		Since $G$ is 2-connected, $G-u$ contains a path $Q$ which has an end $s$ in $C$, an end $t$ in $ \cup_{i=1}^{3}\vts{P_i} $ and no other vertex in these two graphs.
		Without loss of generality, we may assume that $t\in P_1$.
		
		Let $P$ be the $tv$-path  of $P_1$ and let $R$ be the $us$-path of $C$ whose parity is the one of $|\eds{P}|+|\eds{Q}|$. It is easy to check that $P\cup Q\cup R\cup P_2\cup P_3$ is an \oktp of $G$ (with ends $u$ and $v$).
		
	\end{case}
	
	\begin{case} $k=2$.
		Without loss of generality, we may assume that $P_2$ and $P_3$ intersect $C$ at the same vertex $u$ and that $P_1$ meets $C$ at a vertex $s\neq u$.
		Let $Q$ be the odd $su$-path of $C$. Clearly, $Q\cup(\cup_{i=1}^{3}P_i)$ is an \oktp (with ends $u$ and $v$).
	\end{case}
	
	\begin{case} $k=3$.
		Let $Q_1$, $Q_2$ and $Q_3$ be the three paths partitioning the edge-set of $C$ defined by the respective ends of $P_1$, $P_2$ and $P_3$ on $C$. If one of the $Q_i$ is even then, using that $C$ is odd, it is straightforward to check that $C\cup P_1\cup P_2\cup P_3$ contains an \oktp.
		Therefore, we may assume that $Q_1$, $Q_2$ and $Q_3$ are odd. Hence, $C\cup P_1\cup P_2\cup P_3$ is a \toskf.
	\end{case}
	
	In each case we found an \oktp or a \toskf, and this proves the proposition.
\end{dem}

Using an efficient algorithm for finding two vertex-disjoint paths, it is easy to convert this proof into a polynomial-time algorithm which finds an \oktp or a \toskf as stated in the proposition.

The proof of \myref{line-lem-firstEar} uses the following theorem of Frank:

\begin{thm}[Frank \cite{Frank1993}]\label{line-thm-groupFrank}
	Let $G$ be a 2-connected graph. For each edge $e$ of $G$, there exists an optimal ear-decomposition of $G$ whose first ear contains $e$.
	
	Furthermore, such a decomposition can be found in polynomial-time.
\end{thm}



\noindent
The other main ingredient is the following lemma, which may be of independent interest:

\begin{restatable}{lem}{firstEar}\label{line-lem-firstEar}
Each 2-connected non-bipartite graph has an  optimal ear-decomposition whose first ear is an odd circuit.
\end{restatable}

\begin{dem}
	Let $(C,P_1,\ldots,P_k)$ be an optimal ear-decomposition of $G$. If $C$ is odd, then we are done.
	
	Hence, we may assume that $C$ is even. Let $i$ be the smallest integer of $\cro{k}$ such that $C\cup P_1\cup\cdots\cup P_i$ is non-bipartite.
	
	Put $H:=C\cup P_1\cup\cdots\cup P_i$ and let $e\in\eds{P_i}$.
	
	Since $H$ has an ear-decomposition, it is 2-connected. Hence, \myref{line-thm-groupFrank} shows that $H$ has an optimal ear-decomposition $(D,Q_1,\ldots,Q_i)$ whose first ear contains $e$ (the number of ears is indeed $i+1$ as all ear-decompositions of $H$ have the same number of ears).
	
	Clearly, $H-e$ is bipartite. Hence, every circuit of $H$ containing $e$ is odd. In particular, $D$ is odd.
	
	Since $(C,P_1,\ldots,P_k)$ is an optimal ear-decomposition of $G$, the decomposition $(C,P_1,\ldots,P_i)$ must be optimal for $H$.
	
	Hence, the ear-decomposition $(D,Q_1,\ldots,Q_i,P_{i+1},\ldots,P_k)$ is optimal for $G$. This proves the lemma.
\end{dem}

This proof and \myref{line-thm-groupFrank} directly show that \emph{such a decomposition can be found in polynomial-time.}


The last tool is the following easy part of \myref{line-prop-oddInter}:
\begin{prop}\label{line-prop-oddInterPartial}
If a 2-connected graph $G$ has two odd circuits which have at most one common vertex, then $G$ contains an \oktp.
\end{prop}

The \emph{ends} of an \oktp (or an odd theta) are its two vertices of degree 3.  We now prove \myref{line-thm-PhiTwoImpliesK4}.

\begin{dem}[of \myref{line-thm-PhiTwoImpliesK4}]
Clearly, any ear-decomposition of a \toskf which starts with an odd circuit has two odd ears. This shows i)=>ii).

To prove the converse, we may obviously assume that $G$ is non-bipartite.
Suppose that $\phib{G}\geq 2$. We will show a \toskf in $G$.

Since $G$ is 2-connected and non-bipartite, \myref{line-lem-firstEar} shows that $G$ has an optimal ear-decomposition $(C,P_1,\ldots,P_k)$ such that $C$ is odd.

Let $H$ be the graph obtained from $G$ by identifying the vertices of $C$ into a single vertex $c$, keeping the possibly new parallel edges and deleting the loops.

\begin{claim}
$H$ is bipartite.
\end{claim}

Suppose to the contrary that $H$ contains an odd circuit $D$. In $G$, the graph $D$ is either an odd circuit meeting $C$ in at most one vertex or an odd path which has exactly its ends in $C$.

If $D$ is an odd circuit in $G$, \myref{line-prop-oddInterPartial} directly shows an \oktp which contradicts the assumptions on $G$.
Hence, $D$ is an odd path which has exactly its ends in $C$. Therefore, $D\cup C$ is an \oktp: a contradiction.

This ends the proof of Claim 1.

\begin{claim}
$H$ contains an odd theta $T$ containing $c$.
\end{claim}

Since $\phib{G}\geq 2$, there exists $i\in\cro{k}$ such that $P_i$ is odd. Since $G$ is simple and \oktpf, $P_i$ cannot be an edge with both ends in $C$.
Hence, $P_i$ was not deleted as a loop of $H$ and corresponds to a path or a circuit of $H$ with the same length.

As $H$ is bipartite, $P_i$ cannot be a circuit of $H$. Besides, the ends of $P_i$ must clearly belong to the same block $B$ of $H$.
Clearly, the union of the ears of $(C,P_1,\ldots,P_k)$ which are contained in $B$ define a 2-connected subgraph of $B$ for which $P_i$ is an odd ear, and $B$ must contain $c$. Therefore, \myref{line-prop-oddtheta} shows that $B$ contains an odd theta $T$ containing $c$.

This proves Claim 2, and we now show:
\begin{claim}
$c$ is an end of $T$.
\end{claim}

Suppose to the contrary that $c$ is not an end of $T$. Let $u$ and $v$ be the ends of $T$ and  $Q_1$, $Q_2$ and $Q_3$ be the three (odd) $uv$-paths of $T$. Without loss of generality, we may assume that $c$ is not an end of $Q_1$.

First, suppose that $Q_1$ is not a path of $G$. In this case, $Q_1$ corresponds in $G$ to two vertex-disjoint paths $Q_1'$ and $Q_1''$ joining respectively $u$ and $v$ to vertices $s$ and $t$ of $C$. Since $C$ is odd, the two $st$-paths of $C$ have distinct parities. Using these paths, it is straightforward to check that $T\cup C$ always contains an \oktp with ends $u$ and $v$. This contradicts that $G$ is \oktpf.

Hence, we may assume that $Q_1$ remains a path in $G$. Then, $T$ is an odd theta of $G$ which has exactly one vertex $w$ in common with $C$ in $G$.

Since $G$ is 2-connected, $G-w$ contains a path $P$ which joins a vertex $x$ of $C$ to vertex $y$ of $T$ and which has no other vertex in $C\cup T$.

If $y\in \vts{Q_1}$, then (using that $C$ contains $xw$-paths of both parities) it is easy to find an \oktp in $G$ with ends $u$ and $v$, contradicting that $G$ is \oktpf.
Therefore, we may assume without loss of generality that $y\in Q_2$ and that the $uy$-path of $Q_2$ is odd. Again, it is straightforward to build an \oktp of $G$ (with ends $u$ and $y$): a contradiction.

This ends the proof of Claim 3.

Now, let $c'$ be the other end of $T$. The three paths of $T$ in $H$ correspond to three  inner-disjoint odd $\set{c',\vts{C}}$-paths of $G$. Since $G$ is 2-connected and \oktpf, \myref{line-prop-FindingAtosKfour} shows that $G$ contains a \toskf, as required.
\end{dem}

It is straightforward to convert this proof into \emph{a polynomial-time algorithm deciding whether an \oktpf simple graph contains a \toskf and finding such a subdivision if it exists}.

Recall that a graph $G$ is \oktpf if and only if $\beta(G)\leq 1$. Is it true that graphs with $\beta=2$ must contain a \toskf whenever $\overline{\phi}$ is large ? The graphs $H_k$ given in \myref{subsec-Frank} show that the answer is \emph{negative}.
Indeed,  they satisfy $\beta(H_k)=2$ and have an edge whose deletion yields a bipartite graph. Hence they cannot contain a \toskf.

\subsection{Motivation: Cao's thesis}\label{subsec-Cao}

Let $C_5^{+}$ denote the simple graph obtained by adding a single edge to $C_5$.
A graph is \emph{\cfpfree} if it does not contain a totally odd subdivision of $C_5^{+}$. Clearly, \cfpfree graphs are the same as \oktpf simple graphs.

Cao's thesis \cite{Cao1995} contains several results and statements on \cfpfree graphs. Furthermore, it gives a construction procedure for these graphs.

In this section, we first state a corrected version of the result of \cite{Cao1995} on subdivisions of $K_4$ in \cfpfree graphs, which motivated the statement and proof of \myref{line-thm-PhiTwoImpliesK4} (see \myref{sec-SubdivAndCao}).
Besides, we show that some statements on \cfpfree graphs of \cite{Cao1995} and the procedure for their construction are incorrect.

We first recall the definitions of \cite{Cao1995} to keep the same terminology. A graph is \emph{critical non-bipartite} if it is non-bipartite and each pair of odd circuits has at least one common edge. A critical non-bipartite graph is furthermore \emph{elementary} if it has an edge whose deletion yields a bipartite graph.

A graph $H$ is \emph{basic} if it is obtained from a graph $G$ by subdividing each edge of $G$ exactly once (that is, each edge of $G$ is replaced by a path of length 2).
A graph is \emph{critical non-basic} if it is not basic and has an edge whose deletion yields a basic graph.

\begin{figure}[h]
	\centering
	\includegraphics[scale=1]{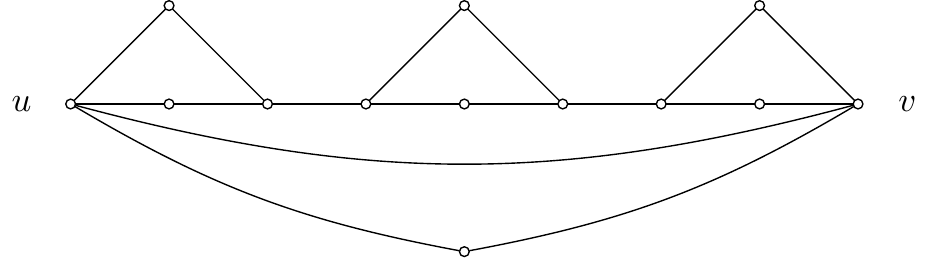}
	\caption[]{an \cfpfree 2-connected elementary critical non-bipartite graph which is not critical non-basic}\label{line-fig-ceCao-elem}
\end{figure}

Clearly, \emph{each critical non-basic graph is \cfpfree and elementary critical non-bipartite}. Lemma 4.5 pg. 70 in \cite{Cao1995} states that the converse also holds: \emph{each 2-connected \cfpfree and elementary critical non-bipartite graph is critical non-basic}. The graph of \myref{line-fig-ceCao-elem} shows that this is false: it is 2-connected, \cfpfree and elementary critical non-bipartite (deleting $uv$ yields a bipartite graph) but it is not critical non-basic.

The following result links totally odd subdivisions of $K_4$ with \cfpfree graphs. In \cite{Cao1995}, it is stated with "critical non-basic" in place of "elementary critical non-bipartite" and the graph of \myref{line-fig-ceCao-elem} shows that it is incorrect as such. Exchanging these two properties corrects the statement:

\begin{thm}[Cao \cite{Cao1995}]
Let $G$ be a non-bipartite graph and $C$ be an odd circuit of $G$.
If $G$ does not contain a \toskf, then for each component $K$ of $G-\eds{C}$: the graph $C\cup K$ is elementary critical non-bipartite.
\end{thm}

Finally, \cite{Cao1995} states a construction procedure for \cfpfree graphs. We observe that it is incorrect. For this purpose, we need only to state a special case of the procedure.

The \emph{sides} of a \toskf are the paths corresponding to the original edges of $K_4$.

Let $F$ be a \toskf. Let $P_1$ and $P_2$ be two vertex-disjoint paths and for each $ i\in\set{1,2} $, let $u_i$ and $v_i$ be the ends of $P_i$.
Let $G$ be a graph obtained by identifying $u_1,v_1,u_2,v_2$ to distinct vertices of $F$ such that for each $ i\in\set{1,2} $: $u_i$ and $v_i$ are identified to vertices which are on sides of $F$ which have a common end $w$, and have even distance to $w$ in $F$.

\begin{figure}[h]
	\centering
	\includegraphics[scale=1.3]{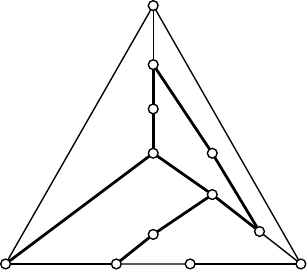}
	\caption{}\label{line-fig-ceCao-K4}
\end{figure}

\cite{Cao1995} states that \emph{each graph obtained in this way is \cfpfree}. The graph of \myref{line-fig-ceCao-K4} shows that this is false: it is obviously built as in the procedure, but the thick edges show an \cfp.

\bibliographystyle{plain}

\bibliography{bibliography}

\begin{thebibliography}{10}

\bibitem{Barahona1994}
F.~Barahona and A.R. Mahjoub.
\newblock Compositions of graphs and polyhedra iii: Graphs with no {$W_4$}
  minor.
\newblock {\em SIAM J. Discrete Math.}, 7(3):372--389, 1994.

\bibitem{Baum1981}
S.~Baum and L.~E.~Jr Trotter.
\newblock Integer rounding for polymatroid and branching optimization problems.
\newblock {\em SIAM Journal on Algebraic and Discrete Methods}, 2:416--425,
  1981.

\bibitem{Benchetrit}
Y.~Benchetrit.
\newblock Integer round-up property for the chromatic number of certain
  h-perfect graphs.
\newblock arXiv:1406.0757, 2014.

\bibitem{Bienstock1991}
D.~Bienstock.
\newblock On the complexity of testing for odd holes and induced odd paths.
\newblock {\em Discrete Mathematics}, 90(1):85 -- 92, 1991.

\bibitem{Bruhn2014}
H.~Bruhn and O.~Schaudt.
\newblock Claw-free t-perfect graphs can be recognised in polynomial time.
\newblock In Jon Lee and Jens Vygen, editors, {\em Integer Programming and
  Combinatorial Optimization}, volume 8494 of {\em Lecture Notes in Computer
  Science}, pages 404--415. Springer International Publishing, 2014.

\bibitem{Cao1995}
D.~Cao.
\newblock {\em Topics in node packing and coloring}.
\newblock PhD thesis, Georgia Institute of Technology, 1995.

\bibitem{Cao1998}
D.~Cao and G.~L. Nemhauser.
\newblock Polyhedral characterizations and perfection of line graphs.
\newblock {\em Discrete Applied Mathematics}, 81:141 -- 154, 1998.

\bibitem{Cheriyan2001}
J.~Cheriyan, A.~Seb{\H{o}}, and Z.~Szigeti.
\newblock Improving on the 1.5-approximation of a smallest 2-edge connected
  spanning subgraph.
\newblock {\em SIAM J. DISCRET. MATH}, 14:170--180, 2001.

\bibitem{Chvatal1975}
V.~Chv{\'a}tal.
\newblock On certain polytopes associated with graphs.
\newblock {\em Journal of Combinatorial Theory, Series B}, 18:138--154, 1975.

\bibitem{Cornuejols2003}
G.~Cornu{\'e}jols, X.~Liu, and K.~Vu{\v{s}}kovi{\'c}.
\newblock A polynomial algorithm for recognizing perfect graphs.
\newblock In {\em 2013 IEEE 54th Annual Symposium on Foundations of Computer
  Science}, pages 20--20. IEEE Computer Society, 2003.

\bibitem{Coullard1996}
C.~Coullard and L.~Hellerstein.
\newblock Independence and port oracles for matroids, with an application to
  computational learning theory.
\newblock {\em Combinatorica}, 16(2):189--208, 1996.

\bibitem{Edmonds1965}
J.~Edmonds.
\newblock Maximum matching and a polyhedron with $0,1$ vertices.
\newblock {\em J. of Res. the Nat. Bureau of Standards}, 69~B:125--130, 1965.

\bibitem{Edmonds1974}
J.~Edmonds and W.~Pulleyblank.
\newblock Facets of 1-matching polyhedra.
\newblock In {\em Hypergraph Seminar of Columbus}, 1974.

\bibitem{Frank1993}
A.~Frank.
\newblock Conservative weightings and ear-decompositions of graphs.
\newblock {\em Combinatorica}, 13(1):65--81, 1993.

\bibitem{Frank2011}
A.~Frank.
\newblock {\em Connections in Combinatorial Optimization}.
\newblock Oxford Lecture Series in Mathematics and its Applications, 2011.

\bibitem{Fulkerson1972}
D.R Fulkerson.
\newblock Anti-blocking polyhedra.
\newblock {\em Journal of Combinatorial Theory, Series B}, 12(1):50 -- 71,
  1972.

\bibitem{Gerards1998}
A.M.H. Gerards and B.~Shepherd.
\newblock The graphs with all subgraphs t-perfect.
\newblock {\em SIAM J. Discrete Math.}, 11(4):524--545, 1998.

\bibitem{Goldberg1973}
M.K. Goldberg.
\newblock On multigraphs of almost maximal chromatic class (in russian).
\newblock {\em Diskret. Analiz}, 23:3--7, 1973.

\bibitem{Grotschel1986}
M.~Gr{\"o}tschel, L.~Lov{\'a}sz, and A.~Schrijver.
\newblock Relaxations of vertex packing.
\newblock {\em Journal of Combinatorial Theory, Series B}, 40(3):330--343,
  1986.

\bibitem{Groetschel1988}
M.~Gr{\"o}tschel, L.~Lov{\'a}sz, and A.~Schrijver.
\newblock {\em Geometric algorithms and combinatorial optimization}.
\newblock Algorithms and combinatorics. Springer-Verlag, Berlin, New York,
  1988.

\bibitem{Huynh2009}
T.~Huynh.
\newblock {\em The linkage-problem for group-labelled graphs}.
\newblock PhD thesis, University of Waterloo, 2009.

\bibitem{Karp1972}
R.~M. Karp.
\newblock Reducibility among combinatorial problems.
\newblock In {\em Complexity of Computer Computations}, pages 85--103, 1972.

\bibitem{Kawarabayashi2010}
K.~Kawarabayashi, Z.~Li, and B.~Reed.
\newblock Recognizing a totally odd $k_4$-subdivision, parity 2-disjoint rooted
  paths and a parity cycle through specified elements.
\newblock In {\em Proceedings of the Twenty-First Annual ACM-SIAM Symposium on
  Discrete Algorithms}, SODA '10, pages 318--328. Society for Industrial and
  Applied Mathematics, 2010.

\bibitem{Kawarabayashi2011}
K.~Kawarabayashi, B.~Reed, and P.~Wollan.
\newblock The graph minor algorithm with parity conditions.
\newblock In {\em Foundations of Computer Science (FOCS), 2011 IEEE 52nd Annual
  Symposium on}, pages 27--36. IEEE, 2011.

\bibitem{Kilakos1996}
K.~Kilakos and F.~B. Shepherd.
\newblock Subdivisions and the chromatic index of {\it r}-graphs.
\newblock {\em Journal of Graph Theory}, 22(3):203--212, 1996.

\bibitem{Lehman1964}
A.~Lehman.
\newblock A solution of the shannon switching game.
\newblock {\em J. Soc. Indust. Appl. Math.}, 12:687--725, 1964.

\bibitem{Lovasz1972}
L.~Lov{\'a}sz.
\newblock Normal hypergraphs and the perfect graph conjecture.
\newblock {\em Discrete Mathematics}, 2(3):253 -- 267, 1972.

\bibitem{Lovasz1972a}
L.~Lov{\'a}sz.
\newblock A note on factor-critical graphs.
\newblock {\em Studia Sci. Math. Hungar.}, 7:279--280, 1972.

\bibitem{Lovasz1985}
L.~Lov{\'a}sz.
\newblock Computing ears and branchings in parallel.
\newblock In {\em Foundations of Computer Science, 1985}, pages 464--467. IEEE,
  1985.

\bibitem{Lovasz2009}
L.~Lov{\'a}sz and M.D. Plummer.
\newblock {\em Matching Theory}.
\newblock AMS Chelsea, 2009.

\bibitem{Miller1986}
G.~L. Miller and V.~Ramachandran.
\newblock Efficient parallel ear decomposition with applications.
\newblock {\em Manuscript, UC Berkeley, MSRI}, 135, 1986.

\bibitem{Oxley1992}
James Oxley.
\newblock {\em Matroid Theory}.
\newblock 1992.

\bibitem{Padberg1973}
M.~W. Padberg.
\newblock On the facial structure of set packing polyhedra.
\newblock {\em Mathematical Programming}, 5(1):199--215, 1973.

\bibitem{Robbins1939}
H.~E. Robbins.
\newblock A theorem on graphs, with an application to a problem of traffic
  control.
\newblock {\em American Mathematical Monthly}, 46:281--283, 1939.

\bibitem{Roussopoulos1973}
N.~D. Roussopoulos.
\newblock A max {m,n} algorithm for determining the graph {H} from its line
  graph {G}.
\newblock {\em Information Processing Letters}, 2(4):108 -- 112, 1973.

\bibitem{Schmidt2013}
J.~M. Schmidt.
\newblock A simple test on 2-vertex-and 2-edge-connectivity.
\newblock {\em Information Processing Letters}, 113(7):241--244, 2013.

\bibitem{Schrijver2003}
A.~Schrijver.
\newblock {\em Combinatorial Optimization - Polyhedra and Efficiency},
  volume~24 of {\em Algorithms and Combinatorics}.
\newblock Springer, 2003.

\bibitem{Seymour1979a}
P.~D. Seymour.
\newblock On multi-colourings of cubic graphs, and conjectures of {F}ulkerson
  and {T}utte.
\newblock {\em Proceedings of the London Mathematical Society},
  s3-38(3):423--460, 1979.

\bibitem{Seymour1979}
P.D. Seymour.
\newblock Some unsolved problems on one-factorization of graphs.
\newblock In J.A. Bondy and U.S.R Murty, editors, {\em Graph theory and related
  topics}, 1979.

\bibitem{Tholey2006}
T.~Tholey.
\newblock Solving the 2-disjoint paths problem in nearly linear time.
\newblock {\em Theory of Computing Systems}, 39(1):51--78, 2006.

\bibitem{Trotter1977}
Jr. Trotter, L.E.
\newblock Line perfect graphs.
\newblock {\em Mathematical Programming}, 12(1):255--259, 1977.

\bibitem{VantHof2012}
P.~Van~'t Hof, M.~Kami{\'n}ski, and D.~Paulusma.
\newblock Finding induced paths of given parity in claw-free graphs.
\newblock {\em Algorithmica}, 62(1-2):537--563, 2012.

\bibitem{Whitney1932}
H.~Whitney.
\newblock Non-separable and planar graphs.
\newblock {\em Transactions of the American Mathematical Society}, 34:339--362,
  1932.

\end{thebibliography}

\end{document}